\documentclass[12pt]{article}%
\usepackage{amsmath}
\usepackage{amsfonts}
\usepackage{amssymb}
\usepackage{graphicx}
\usepackage{amsfonts}
\usepackage{amssymb}%
\makeatletter
\newfont{\footsc}{cmcsc10 at 8truept}
\newfont{\footbf}{cmbx10 at 8truept}
\newfont{\footrm}{cmr10 at 10truept}
\newtheorem{theorem}{Theorem}

\newtheorem{conjecture}[theorem]{Conjecture}
\newtheorem{corollary}[theorem]{Corollary}

\newtheorem{problem}[theorem]{Problem}
\newtheorem{proposition}[theorem]{Proposition}

\newenvironment{proof}[1][Proof]{\noindent{\textbf {#1}  }}  {\hfill$\Box$\bigskip}

\begin{document}

\title{\textbf{The clique number and the smallest }$Q$-\textbf{eigenvalue of graphs}}
\author{Leonardo de Lima\thanks{Production Engineering Department, Federal Center of
Technological Education Celso Suckow da Fonseca, Rio de Janeiro, Brazil,
\textit{email:leonardo.lima@cefet-rj.br} } , Vladimir
Nikiforov\thanks{Department of Mathematical Sciences, University of Memphis,
Memphis TN 38152, USA; \textit{email:vnikifrv@memphis.edu}} , and Carla
Oliveira\thanks{Mathematics Department and Graduate Program of Technology,
National School of Statistics and Federal Center of Technological Education
Celso Suckow da Fonseca, Rio de Janeiro, Brazil,
\textit{email:carla.oliveira@ibge.gov.br}}}
\date{}
\maketitle

\begin{abstract}
Let $q_{\min}(G)$ stand for the smallest eigenvalue of the signless Laplacian
of a graph $G$ of order $n.$ This paper gives some results on the following
extremal problem:

\emph{How large can }$q_{\min}\left(  G\right)  $\emph{ be if }$G$\emph{ is a
graph of order }$n,$\emph{ with no complete subgraph of order }$r+1?$

It is shown that this problem is related to the well-known topic of making
graphs bipartite. Using known classical results, several bounds on $q_{\min}$
are obtained, thus extending previous work of Brandt for regular graphs.

In addition, using graph blowups, a general asymptotic result about the
maximum $q_{\min}$ is established. As a supporting tool, the spectra of the
Laplacian and the signless Laplacian of blowups of graphs are
calculated.\medskip



\end{abstract}

\section{\label{intro}Introduction}

In this paper we study how large can the smallest signless Laplacian
eigenvalue of graphs with bounded clique number be.\medskip

Arguably the most attractive problems in spectral graph theory are the
extremal ones, with general form like:\medskip

\emph{If }$G$\emph{ is a graph of order }$n,$\emph{ with some property
}$\mathcal{P},$\emph{ how large can its }$k$\emph{'th eigenvalue be}%
$?\medskip$

From this general template, by choosing the property $\mathcal{P},$ the type
of graph matrix, and the value $k$, we can obtain an amazing variety of
concrete spectral problems, ranging from trivial to extremely challenging
ones. The study of such extremal questions is crucial to graph theory, for
they provide a sure way to connect the structure of a graph to its
eigenvalues.\medskip

In this vein, we shall introduce a new extremal problem about the smallest
signless Laplacian eigenvalue of graphs with no complete subgraphs of given
order.\medskip

First, recall a few definitions: Given a graph $G,$ write $A$ for the
adjacency matrix of $G$ and let $D$ be the diagonal matrix of the degrees of
$G$. The \emph{Laplacian} $L\left(  G\right)  $ and the \emph{signless
Laplacian} $Q\left(  G\right)  $ of $G$ are defined as $L\left(  G\right)
=D-A$ and $Q\left(  G\right)  =D+A$. We write $\lambda_{1},\ldots,\lambda_{n}$
and $q_{1},\ldots,q_{n}$ for the eigenvalues of $A$ and $Q\left(  G\right)  $
in descending order, and $\mu_{1},\ldots,\mu_{n}$ for the eigenvalues of
$L\left(  G\right)  $ in ascending order. Occasionally we write $q_{\min}$ and
$\lambda_{\min}$ for $q_{n}$ and $\lambda_{n}.$ For more details on the $Q$
matrix, we refer the reader to \cite{Cve10}.\medskip

Here is our new problem:\medskip

\textbf{Problem A }\emph{Let }$n>r\geq2.$\emph{ How large can }$q_{n}\left(
G\right)  $\emph{ be if }$G$\emph{ is graph of order }$n$\emph{ with no
complete subgraph of order }$r+1$\emph{?\medskip}

Note that Problem A is in the spirit of the classical extremal graph theory,
where the analog of Problem A is answered by the Tur\'{a}n theorem. To state
this theorem, let $e\left(  G\right)  $ denote the number of edges of $G,$
write $K_{r}$ for the complete graph of order $r,$ and write $T_{r}(n)$ for
the complete $r$-partite graph of order $n,$ with parts of size $\left\lfloor
n/r\right\rfloor $ or $\left\lceil n/r\right\rceil .$ \medskip

\textbf{Theorem T (Tur\'{a}n, \cite{Tur41})} \emph{If }$n>r\geq2$\emph{ and
}$G$\emph{ is a }$K_{r+1}$\emph{-free} \emph{graph of order }$n,$\emph{ then}
$e\left(  G\right)  <e\left(  T_{r}\left(  n\right)  \right)  ,$ \emph{unless
}$G=T_{r}\left(  n\right)  .\medskip$

Problem A may seem a like of Theorem T, but this similarity is superficial,
for it turns out that Problem A is a much deeper question, entangled with a
notoriously difficult open problem in extremal graph theory. To substantiate
this claim, let us state a theorem, which at first glance seems out of line.

\begin{theorem}
\label{mt}If $G$ is a graph of order $n$, then one has to remove at least
$q_{n}n/4$ edges to make $G$ bipartite.
\end{theorem}

We shall prove Theorem \ref{mt} in Section \ref{pmt}, but note that Brandt
\cite{Bra98} has already proved the same assertion for regular graphs, by a
different method and with a different terminology. However, the general
Theorem \ref{mt} turns out to be much more useful. The reason is that the
topic of making graphs bipartite has been studied for longtime, with several
usable results, which in view of Theorem \ref{mt} directly apply to Problem A.

This line of research was started with the following conjecture of Erd\H{o}s
\cite{Erd71}:

\begin{conjecture}
\label{con1}Every triangle-free graph of order $n$ can be made bipartite by
removing at most $n^{2}/25$ edges.
\end{conjecture}

Defying 46 years of attacks, Conjecture \ref{con1} is still widely open.
Nonetheless, a few nontrivial results are known (see, e.g., \cite{Bra98},
\cite{EFPS88}, \cite{EGS92}, and \cite{Sud07}), which we shall use below for
partial answers to Problem A.

Conjecture \ref{con1} can be extended for $K_{r}$-free graphs; for example, in
\cite{EFPS88} it was conjectured that a $K_{4}$-free graph can be made
bipartite by deleting at most $n^{2}/9$ edges. This conjecture was fully
proved by Sudakov in \cite{Sud07} -- one of the few definite results in this
area. In Section \ref{tmk} we shall use Sudakov's result to get a corollary
about $q_{\min}$ of $K_{4}$-free graphs.

However, the progress with Problem A along this line can go only so far, and
it is unlikely that it can be reduced to a question about making a graph
bipartite. Indeed, Problem A seems to have its own level of difficulty and its
solution may take a while.

To get started, one can simplify Problem A by restating it for regular
graphs$:\medskip$

\textbf{Problem B }\emph{Let }$n>r\geq2.$\emph{ How large can }$q_{n}\left(
G\right)  $\emph{ be if }$G$\emph{ is a regular }$K_{r+1}$\emph{-free graph of
order }$n?$\emph{\medskip}

This step is well justified, for first, the known upper bounds can be
considerably reduced for regular graphs, and second, it is likely that the
extremal graphs in Problem A are regular or close to regular. Hence, Problem B
may provide useful intuition for Problem A.

Moreover, Brandt \cite{Bra98} has already obtained several results for
$q_{\min}$ of regular $K_{r+1}$-free\emph{ }graphs, albeit stated in different
terms. We shall recall some of these results in due course below.

\subsection{The function $f_{r}\left(  n\right)  $ and its asymptotics}

To study Problem A in a systematic way let us define the function
\[
f_{r}\left(  n\right)  :=\max\left\{  q_{n}\left(  G\right)  :\text{ }G\text{
is a graph of order }n\text{ and }G\text{ contains no }K_{r+1}\right\}  .
\]

With $f_{r}\left(  n\right)  $ in hand, we can give a more formal statement of
Problem A:

\begin{problem}
\label{mpro}For any $r\geq2$ and $n>r,$ find or estimate $f_{r}\left(
n\right)  .$
\end{problem}

Note that the introduction of $f_{r}\left(  n\right)  $ does not advance the
solution of Problem A in any concrete way, yet it allows to clearly see and
track the two main lines of attack: on the one hand, obtaining upper bounds on
$f_{r}\left(  n\right)  $ by proofs, and on the other hand, obtaining lower
bounds on $f_{r}\left(  n\right)  $ by constructions. The ultimate goal is to
close the gap between the upper and lower bounds, which, unfortunately, might
take some time.

Before presenting concrete bounds we shall come up with general asymptotics of
$f_{r}\left(  n\right)  .$ For every $r\geq2,$ let us define the real number
$c_{r}$ as
\[
c_{r}:=\sup\left\{  q_{\min}\left(  G\right)  /v\left(  G\right)  :\text{
}G\text{ is a graph with no }K_{r+1}\right\}  .
\]
Because $q_{n}\left(  G\right)  \leq$ $v\left(  G\right)  -2,$ we see that
$c_{r}$ is well defined. Clearly, the definition of $c_{r}$ implies a simple
universal bound for any $K_{r+1}$-free graph $G$ of order $n:$
\[
q_{n}\left(  G\right)  \leq c_{r}n.
\]
What's more, this bound is asymptotically best possible, as given by the next theorem:

\begin{theorem}
\label{tlim}For every $r\geq2,$ the limit
\[
\lim_{n\rightarrow\infty}\frac{1}{n}f_{r}\left(  n\right)
\]
exists and is equal to $c_{r}.$
\end{theorem}

The proof of Theorem \ref{tlim} is given in Section \ref{pbu}. Here we want to
emphasize that this theorem makes the study of $f_{r}\left(  n\right)  $
rather straightforward, as the asymptotic behavior of $f_{r}\left(  n\right)
$ would be determined if we knew the constants $c_{r}$. Unfortunately, at this
stage we do not know any of the constants $c_{r}$ for $r\geq2.$

Another point to make here is that the proof of Theorem \ref{tlim} uses blowup
of graphs. Since the spectra of the signless Laplacian and the Laplacian of
graph blowups have not been studied in the literature, in Section \ref{sbu} we
shall give a few relevant results.

\subsection{Maximal $q_{\min}$ of triangle-free graphs}

In \cite{EFPS88}, Erd\H{o}s, Faudree, Pach, and Spencer have established that
every triangle-free graph of order $n$ can be made bipartite by removing at
most $n^{2}/18+n/2$ edges. Using Theorems \ref{mt} and \ref{tlim}, we
immediately get the inequality $c_{3}<2/9,$ implying the following general
bound, which is the best one known to the authors:

\begin{corollary}
\label{ttu}If $G$ is a triangle-free graph of order $n,$ then $q_{n}\left(
G\right)  <2n/9.$
\end{corollary}

As for lower bounds, Brandt \cite{Bra98} observed that a good lower bound on
the ratio $q_{n}\left(  G\right)  /n$ can be obtained from the Higman-Sims
graph $H_{100},$ introduced by Mesner in 1959 and independently by Higman and
Sims in \cite{HiSi68}. Let us recall that $H_{100}$ is a strongly regular
graph with parameters $\left(  100,22,0,6\right)  .$ Its smallest adjacency
eigenvalue is $-8,$ and since for any $d$-regular graph $H$ we have
\[
q_{\min}\left(  H\right)  =d+\lambda_{\min}\left(  H\right)  ,
\]
for $H_{100}$ we see that%
\[
q_{\min}\left(  H_{100}\right)  =22-8=14.
\]

Clearly, by blowing up $H_{100},$ we obtain a bound for $f_{3}\left(
n\right)  $ for every $n:$

\begin{proposition}
\label{ttl}If $n\geq1,$ then%
\begin{equation}
f_{3}\left(  n\right)  \geq14\left\lfloor n/100\right\rfloor . \label{lob}%
\end{equation}

\end{proposition}

Bound (\ref{lob}) is the best lower bound that we are aware of. The fact that
this bound is based on such complicated graph as the Higman-Sims graph leaves
us clueless as to what $c_{3}$ might be. As Brandt \cite{Bra98} pointed out,
if Conjecture \ref{con1} is true, then we would have $f_{3}\left(  n\right)
\leq0.16n,$ which is quite close to $0.14n.$ On the other hand, if true,
Conjecture \ref{con1} is best possible, while there is no ground to believe
that the inequality $f_{3}\left(  n\right)  \leq0.16n$ is tight.

For completeness, let us mention that for regular triangle-free graphs Brandt
\cite{Bra98} has shown that
\[
q_{n}\left(  G\right)  \leq\left(  3-2\sqrt{2}\right)  n<0.1716n,
\]
which also is not too far from $0.14n.$

\subsection{\label{tmk}Maximal $q_{\min}$ of $K_{r+1}$-free graphs}

Using the Tur\'{a}n graphs $T_{r}\left(  n\right)  $, one can easily see the
following lower bound:

\begin{proposition}
\label{pro3}If $r\geq3,$ then $f_{r}\left(  n\right)  \geq\left(  r-2\right)
\left\lfloor n/r\right\rfloor .$
\end{proposition}

The bound is also valid for $r=2,$ but is meaningless, so triangle-free graphs
need a separate approach.

For $r\geq3,$ we can give only the following approximation, which certainly is
not tight.

\begin{theorem}
\label{tup}If $r\geq3$ and $G$ is a $K_{r+1}$-free graph of order $n,$ then%
\[
q_{n}\left(  G\right)  <\left(  1-\frac{3}{3r-1}\right)  n.
\]

\end{theorem}

Theorem \ref{tup} is proved in Section \ref{pmt}.

For completeness, let us put the known general bounds on $c_{r}$ on one line.

\begin{corollary}
If $r\geq3,$ then%
\[
1-\frac{2}{r}<c_{r}\leq1-\frac{3}{3r-1}.
\]
$\allowbreak$
\end{corollary}

As noted before, Sudakov \cite{Sud07} proved that every $K_{4}$-free graph of
order $n$ can be made bipartite by deleting at most $n^{2}/9$ edges. Combining
this fact with Theorem \ref{mt}, we get the following upper bound:

\begin{theorem}
If $G$ is a $K_{4}$-free graph of order $n,$ then%
\begin{equation}
q_{n}\left(  G\right)  \leq4n/9. \label{k4in}%
\end{equation}

\end{theorem}

This is the best known upper bound on $q_{\min}\left(  G\right)  $ of $K_{4}%
$-free graphs. Note again, that Sudakov's result is best possible, but there
is no evidence that inequality (\ref{k4in}) is tight.

The authors of this note have investigated quite a few small graphs in search
of maximal $q_{\min}$ as a function of the clique number of the graph.
Eventually we believe that the following conjecture might hold:

\begin{conjecture}
\label{con2}Let $r\geq3$ and let $n$ be sufficiently large. If $G$ is a
$K_{r+1}$-free graph of order $n,$ then
\[
q_{n}\left(  G\right)  <q_{n}\left(  T_{r}\left(  n\right)  \right)  ,
\]
unless $G=T_{r}\left(  n\right)  .$
\end{conjecture}

Clearly, if true, Conjecture \ref{con2} is best possible. It is not hard to
see that
\[
\left(  r-2\right)  \left\lfloor n/r\right\rfloor <q_{n}\left(  T_{r}\left(
n\right)  \right)  \leq\left(  1-\frac{2}{r}\right)  n.
\]

Therefore, we have the following weaker form of Conjecture \ref{con2}:

\begin{conjecture}
\label{con3}If $r\geq3,$ then $c_{r}=1-2/r.$
\end{conjecture}

Some credibility to this conjecture is given by the following result of Brandt
\cite{Bra98} about regular graphs:

\begin{theorem}
\label{tB}If $r\geq3$ and $G$ is a $K_{r+1}$-free regular graph of order $n,$
then%
\begin{equation}
q_{n}\left(  G\right)  \leq\left(  5-\frac{4}{r}\left(  \sqrt{r^{2}%
-r}+1\right)  \right)  n. \label{bB}%
\end{equation}

\end{theorem}

Indeed, it is not hard to see that
\[
1-2/r<5-\frac{4}{r}\left(  \sqrt{r^{2}-r}+1\right)  <1-2/r+O\left(
r^{-2}\right)  ,
\]
so if $r$ tends to infinity, bound \ref{bB} approaches the best possible one.

\subsection{\label{sbu}Laplacians and signless Laplacians of graph blowups}

Given a graph $G$ and an integer $t\geq1,$ write $G^{\left(  t\right)  }$ for
the graph obtained by replacing each vertex $u$ of $G$ by a set $V_{u}$ of $t$
independent vertices and every edge $\left\{  u,v\right\}  $ of $G$ by a
complete bipartite graph with parts $V_{u}$ and $V_{v}.$ Usually $G^{\left(
t\right)  }$ is called a \textit{blowup} of $G.$ Blowups of graphs are a very
important tool in the extremal and structural theories of graphs and
hypergraphs, see, e.g., the classical work of Sidorenko \cite{Sid87}.

Although introduced by a combinatorial definition, graph blowups have a clear
algebraic meaning as well: if $A$ is the adjacency matrix of $G,$ then the
adjacency matrix $A\left(  G^{\left(  t\right)  }\right)  $ of $G^{\left(
t\right)  }$ is given by
\[
A\left(  G^{\left(  t\right)  }\right)  =A\otimes J_{t},
\]
where $\otimes$ is the Kronecker product and $J_{t}$ is the all-ones square
matrix of order $t.$ This observation yields the following facts (see, e.g.,
\cite{Nik06b}).

\begin{proposition}
\label{pro1}The\ eigenvalues of $G^{\left(  t\right)  }$ are $t\lambda
_{1}\left(  G\right)  ,\ldots,t\lambda_{n}\left(  G\right)  ,$ together with
$n\left(  t-1\right)  $ additional $0$'s.
\end{proposition}

We write $\overline{G}$ for the complement of a graph $G.$

\begin{proposition}
\label{pro2}The\ eigenvalues of $\overline{G^{\left(  t\right)  }}$ are
$t\lambda_{1}(\overline{G})+t-1,\ldots,t\lambda_{n}(\overline{G})+t-1,$
together with $n\left(  t-1\right)  $ additional $-1$'s.
\end{proposition}

The algebraic meaning of graph blowups make them equally important in spectral
graph theory, see, e.g., \cite{Nik06b} and \cite{Nik15}.

However, the Laplacian and the signless Laplacian of graph blowups are not so
immediately related to the Kronecker product and have not been considered in
the literature as yet. This is unfortunate as there are many spectral problems
about the Laplacian and the signless Laplacian that might benefit from the
blowup construction if theorems similar to Propositions \ref{pro1} and
\ref{pro2} are available.

The goal of this section is to state such basic results about the spectra of
$L(G^{(t)})$, $Q(G^{(t)}),$ and $Q(\overline{G^{(t)}}),$ which, however, turn
out to be more difficult than for the adjacency matrix$.$

Let us also add that the spectrum of $L(\overline{G^{(t)}})$ is obtained
immediately from our results, as $\mu_{n-i}(\overline{G})=n-\mu_{i}\left(
G\right)  ,$ for $i=1,\ldots,n-1,$ (see, e.g., \cite{AM85}).

Here are our three theorems:

\begin{theorem}
\label{thm1} Let $t\geq2.$ If $G$ is a graph of order $n$, with degrees
$d_{1},\ldots,d_{n}$ and Laplacian eigenvalues ${\mu}_{1},\ldots,{\mu}_{n},$
then the Laplacian eigenvalues of $G^{(t)}$ are
\[
t{\mu}_{1},\ldots,t{\mu}_{n},td_{1},\ldots,td_{n},
\]
where each of the eigenvalues $td_{1},\ldots,td_{n}$ has multiplicity $t-1$.
\end{theorem}

The proof of Theorem \ref{thm1} is given in Section \ref{pbu}.

\begin{theorem}
\label{thm2} Let $t\geq2.$ If $G$ is a graph of order $n$, with degrees
$d_{1},\ldots,d_{n},$ and signless Laplacian eigenvalues ${q}_{1},\ldots
,q_{n},$ then the signless Laplacian eigenvalues of $G^{(t)}$ are
\[
tq_{1},\ldots,tq_{n},td_{1},\ldots,td_{n},
\]
where each of the eigenvalues $td_{1},\ldots,td_{n}$ has multiplicity $t-1$.
\end{theorem}

The proof of Theorem \ref{thm1} works with minor changes for Theorem
\ref{thm2}, so we shall omit it.

Here is our final theorem, about the signless Laplacian of the complement of a
blowup. Its proof is also in Section \ref{pbu}.

\begin{theorem}
\label{thm4}Let $t\geq2$ and let $G$ be a graph of order $n$, with degrees
$d_{1},\ldots,d_{n}.$ If the eigenvalues of the signless Laplacian of
$\overline{G}$ are $\overline{q}_{1},\ldots,\overline{q}_{n},$ then the
eigenvalues of the signless Laplacian of $\overline{G^{\left(  t\right)  }}$
are
\[
t\overline{q}_{1}+2(t-1),\ldots,t\overline{q}_{n}+2(t-1),tn-td_{1}%
-2,\ldots,tn-td_{n}-2,
\]
where each of the eigenvalues $tn-td_{1}-2,\ldots,tn-td_{n}-2$ has
multiplicity $t-1$.
\end{theorem}

\section{\label{main}Proofs}

For general terminology and notation on graphs we refer the reader to
\cite{Bol98}.

As usual we write $I_{n}$ and $J_{n}$ for the identity and the all-ones
matrices of order $n.$

\subsection{\label{pmt}Proofs of Theorems \ref{mt} \ and \ \ref{tup}}

\begin{proof}
[\textbf{Proof of Theorem \ref{mt}}]Let $G=G\left(  V,E\right)  $ be a graph
of order $n$ with vertex set $V$ and edge set $E.$ Let $H$ be a bipartite
subgraph of $G$ of maximal number of edges. Write $A$ and $B$ for the vertex
classes of $H$ and let $C=V\backslash\left(  A\cup B\right)  .$ We shall show
that $C$ is either empty or consists of isolated vertices of $G.$ Let us write
$\Gamma\left(  u\right)  $ for the set of neighbors of a vertex $u$ of $G.$

Let $u\in C.$ If $\Gamma\left(  u\right)  \cap A\neq\varnothing,$ add $u$ to
$B$ and all edges joining $u$ to $A$ to $E\left(  H\right)  $. The resulting
graph is bipartite and has more edges than $H,$ contradicting the choice of
$H.$ Hence, for any $u\in C,$ $\Gamma\left(  u\right)  \cap A=\varnothing,$
and by symmetry $\Gamma\left(  u\right)  \cap B=\varnothing.$ Now, if $C$
induces at least one edge $\left\{  u,v\right\}  $ in $G,$ then adding $u$ to
$A$, $v$ to $B,$ and $\left\{  u,v\right\}  $ to $E\left(  H\right)  ,$ we
obtain again a bipartite subgraph of $G$ with more edges than $H,$
contradicting the choice of $H.$ Hence, either $C=\varnothing$ or $C$ consists
of isolated vertices.

Now, take a vector $\mathbf{x}:=\left(  x_{1},\ldots,x_{n}\right)  $ such
that
\[
x_{u}:=\left\{
\begin{array}
[c]{cc}%
1/\sqrt{n} & \text{if }u\in A\text{;}\\
-1/\sqrt{n} & \text{if }u\in V\backslash A.
\end{array}
\right.
\]
Note that if $\left\{  u,v\right\}  $ $\in E\left(  H\right)  ,$ then
$x_{u}+x_{v}=0,$ while if $\left\{  u,v\right\}  \in$ $E\left(  G\right)
\backslash E\left(  H\right)  ,$ then $x_{u}+x_{v}=\pm2/\sqrt{n}.$ Hence,
\[
\sum_{\left\{  u,v\right\}  \in E\left(  G\right)  }\left(  x_{u}%
+x_{v}\right)  ^{2}=\frac{4}{n}\left(  E\left(  G\right)  -E\left(  H\right)
\right)  .
\]
By Rayleigh's principle,
\[
q_{n}\left(  G\right)  =\min\left\{  \left\langle Q\mathbf{x},\mathbf{x}%
\right\rangle :\left\Vert \mathbf{x}\right\Vert =1\right\}  =\min\left\{
\sum_{\left\{  u,v\right\}  \in E\left(  G\right)  }\left(  x_{u}%
+x_{v}\right)  ^{2}:x_{1}^{2}+\cdots+x_{n}^{2}=1\right\}  .
\]
Hence,%
\[
q_{n}\left(  G\right)  n/4\leq E\left(  G\right)  -E\left(  H\right)  ,
\]
proving Theorem \ref{mt}.
\end{proof}

\medskip

\begin{proof}
[\textbf{Proof of Theorem \ref{tup}}]Let $r\geq3$ and let the graph $G$
satisfy the premises of the theorem. Write $\delta$ for the minimum degree of
$G$ and $m$ for the number of its edges. In \cite{LOAN11} it was proved that
if $G$ is $r$-partite, then
\[
q_{\min}\left(  G\right)  \leq\left(  \frac{r-2}{r-1}\right)  \frac{2m}{n},
\]
which, in view of $m\leq\left(  1-1/r\right)  n^{2}/2,$ gives immediately%
\[
q_{\min}\left(  G\right)  \leq\left(  1-\frac{2}{r}\right)  n<\left(
1-\frac{3}{3r-1}\right)  n.
\]
So we shall suppose that the chromatic number of $G$ is at least $r+1.$ On the
other hand, a celebrated theorem of Andr\'{a}sfai, Erd\H{o}s and S\'{o}s
\cite{AES74} shows that if $r\geq2$ and $G$ is a $K_{r+1}$-free graph of order
$n$ and
\[
\delta>\left(  1-\frac{3}{3r-1}\right)  n,
\]
then $G$ is $r$-partite. Since in our case $G$ is not $r$-partite, we conclude
that
\begin{equation}
\delta\leq\left(  1-\frac{3}{3r-1}\right)  n.\label{bAES}%
\end{equation}
Now, recall that in \cite{LOAN11} it was proved that $q_{\min}\left(
G\right)  <\delta.$ Combining this inequality with (\ref{bAES}), completes the
proof of Theorem \ref{tup}.
\end{proof}

\subsection{\label{pbu}Proof of Theorems \ref{thm1}, \ref{thm4}, and
\ref{tlim}}

\begin{proof}
[\textbf{Proof of Theorem \ref{thm1}}]Let $G$ be a graph satisfying the
hypothesis of the theorem, and let $A=\left[  a_{i,j}\right]  $ be its
adjacency matrix, $D$ be the diagonal matrix of its degrees and $L$ be its
Laplacian, i.e. $L=D-A$. With appropriate labeling, the Laplacian matrix of
$G^{(t)}$ can be written as a $t\times t$ block matrix%
\[
B=\left[
\begin{array}
[c]{cccc}%
L+(t-1)D & -A & \ldots & -A\\
-A & L+(t-1)D & \ldots & -A\\
\vdots & \vdots & \ddots & \vdots\\
-A & -A & \ldots & L+(t-1)D
\end{array}
\right]  .
\]
Now, let $\mathbf{x}_{1},\ldots,\mathbf{x}_{n}$ be pairwise orthogonal
eigenvectors to ${\mu}_{1},\ldots,{\mu}_{n}.$ For convenience we represent
$\mathbf{x}_{1},\ldots,\mathbf{x}_{n}$ as column vectors. For each
$i\in\left[  n\right]  ,$ define a column vector $\mathbf{y}_{i}$ of length
$tn$ as
\[
\mathbf{y}_{i}:=\left[
\begin{array}
[c]{c}%
\mathbf{x}_{i}\\
\vdots\\
\mathbf{x}_{i}%
\end{array}
\right]  ,
\]
and note that%
\[
B\mathbf{y}_{i}=\left[
\begin{array}
[c]{c}%
L\mathbf{x}_{i}+(t-1)D\mathbf{x}_{i}-(t-1)A\mathbf{x}_{i}\\
\vdots\\
L\mathbf{x}_{i}+(t-1)D\mathbf{x}_{i}-(t-1)A\mathbf{x}_{i}%
\end{array}
\right]  =\left[
\begin{array}
[c]{c}%
tL\mathbf{x}_{i}\\
\vdots\\
tL\mathbf{x}_{i}%
\end{array}
\right]  =\left[
\begin{array}
[c]{c}%
t{\mu}_{i}\mathbf{x}_{i}\\
\vdots\\
t{\mu}_{i}\mathbf{x}_{i}%
\end{array}
\right]  =t\mu_{i}\mathbf{y}_{i}.
\]
Hence $t\mu_{i}$ is an eigenvalue of $L(G^{(t)})$ with eigenvector
$\mathbf{y}_{i}.$ Clearly $\mathbf{y}_{1},\ldots,\mathbf{y}_{n}$ are pairwise
orthogonal, as $\mathbf{x}_{1},\ldots,\mathbf{x}_{n}$ are pairwise orthogonal.

To find the remaining $nt-n$ eigenvalues of $L(G^{(t)}),$ fix some
$s\in\left[  n\right]  $ and write $\mathbf{e}_{j}$ for the column vector of
length $tn$ having $1$ at position $j$ and zeros elsewhere. For $k=1,\ldots
,t-1$, define a column vector $\mathbf{z}_{k}^{s}=\mathbf{e}_{s}%
-\mathbf{e}_{kn+s}$. Note that $B\mathbf{e}_{s}$ and $B\mathbf{e}_{kn+s}$ are
just the $s$'th and the $\left(  kn+s\right)  $'th columns of $B,$ which
coincide everywhere, but at the $s$'th and the $\left(  kn+s\right)  $'th
entries. Thus, we see that
\[
B\mathbf{z}_{k}^{s}=B\mathbf{e}_{s}-B\mathbf{e}_{kn+s}=td_{s}\left(
\mathbf{e}_{s}-\mathbf{e}_{kn+s}\right)  .
\]
Therefore, $\mathbf{z}_{1}^{s},\ldots,\mathbf{z}_{t-1}^{s}$ are eigenvectors
to the eigenvalue $td_{s}.$ We shall show that $\mathbf{z}_{1}^{s}%
,\ldots,\mathbf{z}_{t-1}^{s}$ are linearly independent$.$ Suppose that
$\sum_{j=1}^{t-1}c_{j}\mathbf{z}_{j}^{s}=0$. Hence%
\[
\left(  \sum_{k=1}^{t-1}c_{k}\right)  \mathbf{e}_{s}=\sum_{k=2}^{t-1}%
c_{k}\mathbf{e}_{kn+s},
\]
which is only possible if $c_{i}=0$ for $i=1,\ldots,t-1.$ Hence $td_{s}$ is an
eigenvalue of multiplicity at least $t-1.$ Furthermore, it is easily seen that
the eigenspaces $Span\left\{  \mathbf{z}_{1}^{p},\ldots,\mathbf{z}_{t-1}%
^{p}\right\}  $ and $Span\left\{  \mathbf{z}_{1}^{q},\ldots,\mathbf{z}%
_{t-1}^{q}\right\}  $ corresponding to $td_{p}$ and to $td_{q}$ are
orthogonal. Finally, for any $i\in\left[  n\right]  $, $s\in\left[  n\right]
$ and $k\in\left[  t-1\right]  $ the vectors $\mathbf{y}_{i}$ and
$\mathbf{z}_{k}^{s}$ are orthogonal, for the $s$'th and the $\left(
kn+s\right)  $'th entries of $\mathbf{y}_{i}$ are the same. Hence, the
eigenspaces corresponding to $t{\mu}_{1},\ldots,t{\mu}_{n},td_{1}%
,\ldots,td_{n}$ are orthogonal and each of the eigenvalues $td_{1}%
,\ldots,td_{n}$ has multiplicity $t-1$. Theorem \ref{thm1} is proved.
\end{proof}

\bigskip

\begin{proof}
[\textbf{Proof of Theorem \ref{thm4}}]Let $G$ be a graph satisfying the
hypothesis of the theorem, and let $\overline{A}$ be the adjacency matrix of
$\overline{G}$, let $\overline{D}$ be the diagonal matrix of the degrees of
$\overline{G},$ and let $\overline{Q}$ be its signless Laplacian, i.e.,
$\overline{Q}=\overline{D}+\overline{A}$. With appropriate labeling, the
signless Laplacian of $\overline{G^{\left(  t\right)  }}$ can be written as a
$t\times t$ block matrix%
\[
B=\left[
\begin{array}
[c]{cccc}%
\overline{Q}+(t-1)(\overline{D}+I_{n}) & \overline{A}+I_{n} & \ldots &
\overline{A}+I_{n}\\
\overline{A}+I_{n} & \overline{Q}+(t-1)(\overline{D}+I_{n}) & \ldots &
\overline{A}+I_{n}\\
\vdots & \vdots & \ddots & \vdots\\
\overline{A}+I_{n} & \overline{A}+I_{n} & \ldots & \overline{Q}%
+(t-1)(\overline{D}+I_{n})
\end{array}
\right]  .
\]
Now, let $\mathbf{x}_{1},\ldots,\mathbf{x}_{n}$ be pairwise orthogonal
eigenvectors to $\overline{{q}}_{1},\ldots,\overline{q}_{n}.$ For convenience
we represent $\mathbf{x}_{1},\ldots,\mathbf{x}_{n}$ as column vectors. For
each $i\in\left[  n\right]  ,$ define a column vector $\mathbf{y}_{i}$ of
length $tn$ as
\[
\mathbf{y}_{i}:=\left[
\begin{array}
[c]{c}%
\mathbf{x}_{i}\\
\vdots\\
\mathbf{x}_{i}%
\end{array}
\right]  ,
\]
and note that%
\begin{align*}
B\mathbf{y}_{i} &  =\left[
\begin{array}
[c]{c}%
\overline{Q}\mathbf{x}_{i}+(t-1)(I_{n}+\overline{D})\mathbf{x}_{i}%
+(t-1)(\overline{A}+I_{n})\mathbf{x}_{i}\\
\vdots\\
\overline{Q}\mathbf{x}_{i}+(t-1)(I_{n}+\overline{D})\mathbf{x}_{i}%
+(t-1)(\overline{A}+I_{n})\mathbf{x}_{i}%
\end{array}
\right]  \\
&  =\left[
\begin{array}
[c]{c}%
{\overline{q}}_{i}\mathbf{x}_{i}+2(t-1)\mathbf{x}_{i}+(t-1)(\overline
{D}+\overline{A})\mathbf{x}_{i}\\
\vdots\\
{\overline{q}}_{i}\mathbf{x}_{i}+2(t-1)\mathbf{x}_{i}+(t-1)(\overline
{D}+\overline{A})\mathbf{x}_{i}%
\end{array}
\right]  \\
&  =\left[
\begin{array}
[c]{c}%
(t{\overline{q}}_{i}+2(t-1))\mathbf{x}_{i}\\
\vdots\\
(t{\overline{q}}_{i}+2(t-1))\mathbf{x}_{i}%
\end{array}
\right]  \\
&  =(t{\overline{q}}_{i}+2(t-1))\mathbf{y}_{i}.
\end{align*}
Hence $t{\overline{q}}_{i}+2(t-1)$ is an eigenvalue of $\overline{G^{\left(
t\right)  }}$ with eigenvector $\mathbf{y}_{i}.$ Clearly $\mathbf{y}%
_{1},\ldots,\mathbf{y}_{n}$ are pairwise orthogonal, as $\mathbf{x}_{1}%
,\ldots,\mathbf{x}_{n}$ are pairwise orthogonal.

To find the remaining $nt-n$ eigenvalues of $Q(\overline{G^{(t)}}),$ let us
note that $B$ can be written as
\[
B=\left[
\begin{array}
[c]{ccc}%
J_{n}-A+(tn-2)I_{n}-tD & \ldots & J_{n}-A\\
\vdots & \ddots & \vdots\\
J_{n}-A & \ldots & J_{n}-A+(tn-2)I_{n}-tD
\end{array}
\right]  .
\]
Now, fix some $s\in\left[  n\right]  $ and write $\mathbf{e}_{j}$ for the
column vector of length $tn$ having $1$ at position $j$ and zeros elsewhere.
For $k=1,\ldots,t-1$, define a column vector $\mathbf{z}_{k}^{s}%
=\mathbf{e}_{s}-\mathbf{e}_{kn+s}$. Note that $B\mathbf{e}_{s}$ and
$B\mathbf{e}_{kn+s}$ are just the $s$'th and the $\left(  kn+s\right)  $'th
columns of $B,$ which coincide everywhere, but at the $s$'th and the $\left(
kn+s\right)  $'th entries. Thus, we see that
\[
B\mathbf{z}_{k}^{s}=B\mathbf{e}_{s}-B\mathbf{e}_{kn+s}=(tn-2-td_{s})\left(
\mathbf{e}_{s}-\mathbf{e}_{kn+s}\right)  .
\]
Therefore, $\mathbf{z}_{1}^{s},\ldots,\mathbf{z}_{t-1}^{s}$ are eigenvectors
to the eigenvalue $tn-2-td_{s}.$ We shall show that $\mathbf{z}_{1}^{s}%
,\ldots,\mathbf{z}_{t-1}^{s}$ are linearly independent$.$ Suppose that
$\sum_{j=1}^{t-1}c_{j}\mathbf{z}_{j}^{s}=0$. Hence%
\[
\left(  \sum_{k=1}^{t-1}c_{k}\right)  \mathbf{e}_{s}=\sum_{k=2}^{t-1}%
c_{k}\mathbf{e}_{kn+s},
\]
which is only possible if $c_{i}=0$ for $i=1,\ldots,t-1.$ Hence $tn-2-td_{s}$
is an eigenvalue of multiplicity at least $t-1.$ Furthermore, it is easily
seen that the eigenspaces $Span\left\{  \mathbf{z}_{1}^{p},\ldots
,\mathbf{z}_{t-1}^{p}\right\}  $ and $Span\left\{  \mathbf{z}_{1}^{q}%
,\ldots,\mathbf{z}_{t-1}^{q}\right\}  $ corresponding to $tn-2-td_{p}$ and to
$tn-2-td_{q}$ are orthogonal. Finally, for any $i\in\left[  n\right]  $,
$s\in\left[  n\right]  $ and $k\in\left[  t-1\right]  ,$ the vectors
$\mathbf{y}_{i}$ and $\mathbf{z}_{k}^{s}$ are orthogonal, for the $s$'th and
the $\left(  kn+s\right)  $'th entries of $\mathbf{y}_{i}$ are the same.
Hence, the eigenspaces corresponding to
\[
t\overline{q}_{1}+2(t-1),\ldots,t\overline{q}_{n}+2(t-1),tn-2-td_{1}%
,\ldots,tn-2-td_{n}%
\]
are orthogonal and each of the eigenvalues $tn-2-td_{1},\ldots,tn-2-td_{n}$
has multiplicity $t-1$. Theorem \ref{thm4} is proved.
\end{proof}

\bigskip

\begin{proof}
[\textbf{Proof of Theorem \ref{tlim}}]Note that if $G$ is a graph and
$t\geq2,$ then%
\begin{equation}
q_{\min}\left(  G^{\left(  t\right)  }\right)  =tq_{\min}\left(  G\right)  .
\label{inq}%
\end{equation}
Indeed write $n$ for the order of $G$ and let $q_{1},\ldots,q_{n}$ and
$d_{1},\ldots,d_{n}$ be the signless Laplacian eigenvalues and the degrees of
$G.$

The eigenvalues of $Q\left(  G^{\left(  t\right)  }\right)  $ are
$tq_{1},\ldots,tq_{n},td_{1},\ldots,td_{n}.$ Since $q_{n}\leq\delta\left(
G\right)  ,$ we see that
\[
q_{n}=\min\left\{  q_{1},\ldots,q_{n},d_{1},\ldots,d_{n}\right\}  ,
\]
and so
\[
tq_{n}=\min\left\{  tq_{1},\ldots,tq_{n},td_{1},\ldots,td_{n}\right\}
=q_{\min}\left(  G^{\left(  t\right)  }\right)  .
\]

Note that the definition of $c_{r}$ implies that for every $n,$
\[
\frac{1}{n}f_{r}\left(  n\right)  \leq c_{r}.
\]
So to prove the assertion we need to show that for all sufficiently large $n,$
one has%
\[
\frac{1}{k}f_{r}\left(  k\right)  \geq c_{r}-\varepsilon.
\]
Choose a graph $G,$ say of order $n,$ such that
\[
\frac{1}{n}q_{n}\left(  G\right)  \geq c_{r}-\varepsilon/2.
\]
We shall show that if $k>2nc_{r}/\varepsilon,$ then
\[
\frac{1}{k}f_{r}\left(  k\right)  \geq c_{r}-\varepsilon.
\]
Indeed, this is obvious if $k$ is a multiple of $n,$ for if $k=nt,$ then
\[
\frac{1}{tn}f_{r}\left(  tn\right)  \geq\frac{1}{tn}q_{n}\left(  G^{\left(
t\right)  }\right)  =\frac{1}{n}q_{n}\left(  G\right)  \geq c_{r}%
-\varepsilon/2.
\]
Now, let $nt<k<n\left(  t+1\right)  ,$ and let $H$ be the union of $G^{\left(
t\right)  }$ and $\left(  k-nt\right)  $ isolated vertices. Clearly
\[
q_{\min}\left(  H\right)  =q_{\min}\left(  G^{\left(  t\right)  }\right)  ,
\]
and so,%
\[
\frac{1}{k}q_{\min}\left(  H\right)  >\frac{1}{n\left(  t+1\right)  }q_{\min
}\left(  G^{\left(  t\right)  }\right)  =\frac{t}{n\left(  t+1\right)  }%
q_{n}\left(  G\right)  \geq\frac{t}{t+1}\left(  c_{r}-\varepsilon/2\right)  .
\]
Obviously, if $k>2nc_{r}/\varepsilon,$ then $t+1>2c_{r}/\varepsilon$ and some
simple algebra shows that
\[
\frac{t}{t+1}\left(  c_{r}-\varepsilon/2\right)  \geq c_{r}-\varepsilon.
\]
Hence, if $k>2c_{r}n/\varepsilon,$ then
\[
\frac{1}{k}f_{r}\left(  k\right)  \geq c_{r}-\varepsilon,
\]
completing the proof of Theorem \ref{tlim}.
\end{proof}

\section{Concluding remarks}

(A) One may wonder if problems similar to Problem A arise if $q_{n}\left(
G\right)  $ is replaced by $\mu_{2}\left(  G\right)  ,$ also known as the
\emph{algebraic connectivity} of $G,$ e.g.:\textbf{ \ \medskip}

\emph{Let }$n>r\geq2.$\emph{ How large can }$\mu_{2}\left(  G\right)  $\emph{
be if }$G$\emph{ is graph of order }$n$\emph{ with no complete subgraph of
order }$r+1$\emph{?\medskip}

However, this question is almost trivial, as for the graph $T_{r}\left(
n\right)  $ we have
\[
\mu_{2}\left(  T_{r}\left(  n\right)  \right)  =\delta\left(  T_{r}\left(
n\right)  \right)  =n-\left\lceil n/r\right\rceil ,
\]
and so the answer is $n-\left\lceil n/r\right\rceil ,$ by the Tur\'{a}n
theorem.\medskip

(B) Improving Theorem \ref{mt} is an independent challenge, as the constant
$1/4$ certainly can be increased. So we raise the following problem:

\begin{problem}
What is the supremum of the set of all numbers $C$ such that there exists a
graph $G$ that cannot be made bipartite by deleting fewer than $Cq_{\min
}\left(  G\right)  v\left(  G\right)  $ edges?
\end{problem}

(C) Finally, if Conjecture \ref{con2} is true, then it would open a very
interesting field of investigation along the lines of the classical extremal
graph theory. One peculiarity of $q_{\min}\left(  G\right)  $ is that it
depends more on the distribution of the edges of a graph than on their number,
so it may become a useful tool in extremal graph theory.\bigskip

\textbf{Acknowledgement. }This work was started during the second author's
visit at the Universidade Federal do Rio de Janeiro and Centro Federal de
Educa\c{c}\~{a}o Tecnol\'{o}gica Celso Suckow da Fonseca in the Spring of
2014. He is grateful to his hosts for all the wonderful experience during his stay.

The research of the first author was supported by CNPq Grant 305867/2012--1
and FAPERJ 102.218/2013

The research of the third author was supported by CNPq Grant
305454/2012--9.

\end{document}